\documentclass[a4paper,12pt]{article}
\usepackage{latexsym}
\usepackage{amsmath}
\usepackage{amssymb}
\usepackage{enumerate}

\usepackage{theorem}
\newtheorem{theorem}{Theorem}[section]

\newtheorem{lemma}[theorem]{Lemma}
\newtheorem{corollary}[theorem]{Corollary}

\theorembodyfont{\rmfamily}

\newtheorem{definition}[theorem]{Definition}

\makeatletter

\@addtoreset{equation}{section}
\makeatother

\newcommand{\qedd}{\hfill \Box}
\newcommand{\ve}{\varepsilon}
\newcommand{\del}{\partial}
\newcommand{\lra}{\longrightarrow}

\newcommand{\R}{\ensuremath{\mathbb{R}}}

\newcommand{\cC}{\ensuremath{\mathcal{C}}}
\newcommand{\cF}{\ensuremath{\mathcal{F}}}
\newcommand{\cH}{\ensuremath{\mathcal{H}}}
\newcommand{\bC}{\ensuremath{\mathbf{C}}}
\newcommand{\bL}{\ensuremath{\mathbf{L}}}

\def\vol{\mathop{\mathrm{vol}}\nolimits}
\def\Ric{\mathop{\mathrm{Ric}}\nolimits}
\def\trace{\mathop{\mathrm{trace}}\nolimits}
\def\CD{\mathop{\mathsf{CD}}\nolimits}

\newcommand{\rev}[1]{\overleftarrow{#1}}

\title{Weighted Ricci curvature estimates\\ for Hilbert and Funk geometries}
\author{Shin-ichi Ohta\thanks{Department of Mathematics, Kyoto University,
Kyoto 606-8502, Japan ({\sf sohta@math.kyoto-u.ac.jp});
Supported in part by the Grant-in-Aid for Young Scientists (B) 23740048.}}
\date{}
\begin{document}

\maketitle

\begin{abstract}
We consider Hilbert and Funk geometries on a strongly convex domain in the Euclidean space.
We show that, with respect to the Lebesgue measure on the domain,
Hilbert (resp.\ Funk) metric has the bounded (resp.\ constant negative) weighted Ricci curvature.
As one of corollaries, these metric measure spaces satisfy the curvature-dimension condition
in the sense of Lott, Sturm and Villani.
\end{abstract}

\section{Introduction}

Hilbert~\cite{Hi} introduced the distance function $d_{\cH}$ on
a bounded convex domain $D \subset \R^n$, related to his fourth problem.
Given distinct points $x,y \in D$, denoting by $x'=x+s(y-x)$ and $y'=x+t(y-x)$
the intersections of the boundary $\del D$ and the line passing through $x$ and $y$ with $s<0<t$
(see Figure), Hilbert's distance $d_{\cH}$ is given by
\[ d_{\cH}(x,y)=\frac{1}{2}\log\bigg( \frac{|x'-y| \cdot |x-y'|}{|x'-x| \cdot |y-y'|} \bigg), \]
where $|\cdot|$ stands for the Euclidean norm.
This is indeed a distance function on $D$, and satisfies the interesting property
that line segments between any points are minimizing.
In the particular case where $D$ is the unit ball, $(D,d_{\cH})$ coincides with
the Klein model of the hyperbolic space.
The structure of $(D,d_{\cH})$ has been investigated from geometric and dynamical aspects
(see, for example, \cite{Eg}, \cite{Be}, \cite{CV}).
For instance, $(D,d_{\cH})$ is known to be Gromov hyperbolic
under mild smoothness and convexity assumptions on $D$.

Funk~\cite{Fu} introduced a non-symmetrization of $d_{\cH}$, namely
\[ d_{\cF}(x,y)=\log\bigg( \frac{|x-y'|}{|y-y'|} \bigg). \]
Note that $d_{\cF}(x,y) \neq d_{\cF}(y,x)$, while the triangle inequality
$d_{\cF}(x,z) \le d_{\cF}(x,y)+d_{\cF}(y,z)$ still holds.
Clearly we have $2d_{\cH}(x,y)=d_{\cF}(x,y)+d_{\cF}(y,x)$,
and line segments are minimizing also with respect to Funk's distance.

\begin{center}
\begin{picture}(350,150)
\put(170,5){Figure}

\qbezier(250,80)(245,50)(210,40)
\qbezier(210,40)(180,35)(140,50)
\qbezier(140,50)(110,60)(110,80)
\qbezier(110,80)(110,115)(130,130)
\qbezier(130,130)(180,165)(230,120)
\qbezier(230,120)(247,105)(250,80)
\put(171.5,41){\line(-1,2){43.5}}

\thicklines
\put(140,50){\line(1,1){79}}
\put(158,68){\vector(-1,2){10}}

\put(130,37){$x'$}
\put(163,63){$x$}
\put(190,90){$y$}
\put(220,137){$y'$}
\put(135,80){$v$}
\put(170,30){$a$}
\put(120,135){$b$}

\put(157,67){\rule{2pt}{2pt}}
\put(187,97){\rule{2pt}{2pt}}

\end{picture}
\end{center}

If $\del D$ is smooth and $D$ is strongly convex
(in other words, $\del D$ is positively curved; see Definition~\ref{df:Fstr}),
then $d_{\cH}$ and $d_{\cF}$ are realized by the smooth Finsler structures
\begin{equation}\label{eq:FHFF}
F_{\cH}(x,v)=\frac{|v|}{2} \bigg\{ \frac{1}{|x-a|}+\frac{1}{|x-b|} \bigg\},
 \quad F_{\cF}(x,v)=\frac{|v|}{|x-b|} \quad \text{for}\ v \in T_xD=\R^n,
\end{equation}
respectively (cf.\ \cite[\S 2.3]{Shspr}), where $a=x+sv$ and $b=x+tv$ denote the intersections
of $\del D$ and the line passing through $x$ in the direction $v$ with $s<0<t$ (see Figure).
Note that $2F_{\cH}(x,v)=F_{\cF}(x,v)+F_{\cF}(x,-v)$.
A remarkable feature of these metrics is that they have the constant negative flag curvature
$-1$ and $-1/4$, respectively (cf.\ \cite[Theorem~1]{Ok}, \cite[Theorem~12.2.11]{Shspr}),
provided that $n \ge 2$ as a matter of course.
The flag curvature is a generalization of the sectional curvature in Riemannian geometry,
so that it is natural that $(D,d_{\cH})$ and $(D,d_{\cF})$ enjoy properties of negatively curved spaces.

Recently, the theory of the \emph{weighted Ricci curvature} (see Definition~\ref{df:wRic})
for Finsler manifolds equipped with arbitrary measures has been developed
in connection with optimal transport theory.
It turned out that the weighted Ricci curvature is a natural quantity and
quite useful in the study of geometry and analysis on Finsler manifolds
(see \cite{Oint}, \cite{Ospl}, \cite{OShf}, \cite{OSbw}).
The aim of this article is to show that the weighted Ricci curvature for Hilbert and Funk geometries
admits uniform bounds with respect to the Lebesgue measure $m_{\bL}$ restricted on $D$.

\begin{theorem}[Funk case]\label{th:Funk}
Let $D \subset \R^n$ with $n \ge 2$ be a strongly convex domain such that $\del D$ is smooth.
Then $(D,F_{\cF},m_{\bL})$ has the constant negative weighted Ricci curvature as,
for any unit vector $v \in TD$,
\[ \Ric_{\infty}(v) =-\frac{n-1}{4}, \qquad \Ric_N(v) =-\frac{n-1}{4}-\frac{(n+1)^2}{4(N-n)}
 \quad \text{for}\ N \in (n,\infty). \]
\end{theorem}

\begin{theorem}[Hilbert case]\label{th:Hilb}
Let $D \subset \R^n$ with $n \ge 2$ be a strongly convex domain such that $\del D$ is smooth.
Then the weighted Ricci curvature of $(D,F_{\cH},m_{\bL})$ is bounded as,
for any unit vector $v \in TD$,
\[ \Ric_{\infty}(v) \in \big(\! -(n-1),2 \big], \quad
 \Ric_N(v) \in \bigg(\!\! -(n-1)-\frac{(n+1)^2}{N-n},2 \bigg]
 \ \text{for}\ N \in (n,\infty). \]
\end{theorem}

We stress that our estimates are independent of the choice of the domain $D$.
There are several applications (Corollaries~\ref{cr:CD}, \ref{cr:BW})
via the theory of the weighted Ricci curvature.

The article is organized as follows.
After preliminaries for Finsler geometry and the weighted Ricci curvature,
we prove Theorems~\ref{th:Funk}, \ref{th:Hilb} in Sections~\ref{sc:Funk}, \ref{sc:Hilb}, respectively.
We finally discuss applications and remarks in Section~\ref{sc:appl}.

\section{Preliminaries}\label{sc:prel}

We very briefly review the necessary notions in Finsler geometry,
we refer to \cite{BCS}, \cite{Shlec} and \cite{Shspr} for further reading.
Let $M$ be a connected, $n$-dimensional $\cC^{\infty}$-manifold
without boundary such that $n \ge 2$.
Given a local coordinate $(x^i)_{i=1}^n$ on an open set $\Omega \subset M$,
we always use the coordinate $(x^i,v^j)_{i,j=1}^n$ of $T\Omega$ such that
\[ v=\sum_{j=1}^n v^j \frac{\del}{\del x^j}\Big|_x \in T_xM
 \qquad \text{for}\ x \in \Omega. \]

\begin{definition}[Finsler structures]\label{df:Fstr}
A nonnegative function $F:TM \lra [0,\infty)$ is called
a \emph{$\cC^{\infty}$-Finsler structure} of $M$ if the following three conditions hold.
\begin{enumerate}[(1)]
\item(Regularity)
$F$ is $\cC^{\infty}$ on $TM \setminus 0$, where $0$ stands for the zero section.

\item(\emph{Positive $1$-homogeneity})
It holds $F(cv)=cF(v)$ for all $v \in TM$ and $c>0$.

\item(\emph{Strong convexity})
The $n \times n$ matrix
\begin{equation}\label{eq:gij}
\big( g_{ij}(v) \big)_{i,j=1}^n :=
 \bigg( \frac{1}{2}\frac{\del^2 (F^2)}{\del v^i \del v^j}(v) \bigg)_{i,j=1}^n
\end{equation}
is positive-definite for all $v \in TM \setminus 0$.
\end{enumerate}
\end{definition}

For $x,y \in M$, we can define the \emph{distance} from $x$ to $y$ in a natural way by
\[ d(x,y):=\inf_{\eta} \int_0^1 F\big( \dot{\eta}(t) \big) \,dt, \]
where the infimum is taken over all $\cC^1$-curves $\eta:[0,1] \lra M$
with $\eta(0)=x$ and $\eta(1)=y$.
We remark that this distance can be \emph{nonsymmetric} (namely $d(y,x) \neq d(x,y)$),
since $F$ is only positively homogeneous.
A $\cC^{\infty}$-curve $\eta$ on $M$ is called a \emph{geodesic}
if it is locally minimizing and has a constant speed (i.e., $F(\dot{\eta})$ is constant).

Given $v \in T_xM$, if there is a geodesic $\eta:[0,1] \lra M$ with $\dot{\eta}(0)=v$,
then we define the \emph{exponential map} by $\exp_x(v):=\eta(1)$.
We say that $(M,F)$ is \emph{forward complete} if the exponential map is defined on whole $TM$.
If the \emph{reverse} Finsler manifold $(M,\rev{F})$ with $\rev{F}(v):=F(-v)$ is forward complete, 
then $(M,F)$ is said to be \emph{backward complete}.
We remark that $(D,F_{\cH})$ is both forward and backward complete (they are indeed equivalent
since $\rev{F_{\cH}}=F_{\cH}$), while $(D,F_{\cF})$ is only forward complete.

For each $v \in T_xM \setminus 0$, the positive-definite matrix
$(g_{ij}(v))_{i,j=1}^n$ in \eqref{eq:gij} induces the Riemannian structure $g_v$ of $T_xM$ as
\begin{equation}\label{eq:gv}
g_v\bigg( \sum_{i=1}^n a_i \frac{\del}{\del x^i}\Big|_x,
 \sum_{j=1}^n b_j \frac{\del}{\del x^j}\Big|_x \bigg)
 := \sum_{i,j=1}^n a_i b_j g_{ij}(v).
\end{equation}
Note that $g_{cv}=g_v$ for $c>0$.
This inner product is regarded as the best Riemannian approximation of $F|_{T_xM}$ in the direction $v$,
in the sense that the unit sphere of $g_v$ is tangent to that of $F|_{T_xM}$ at $v/F(v)$ up to the second order.
In particular, we have $g_v(v,v)=F(v)^2$.

The \emph{Ricci curvature} (as the trace of the \emph{flag curvature}) for a Finsler manifold
is defined by using the Chern connection.
Instead of giving the precise definition in coordinates,
we explain a useful interpretation due to Shen (see \cite[\S 6.2]{Shlec}, \cite[Lemma~2.4]{Shcdv}).
Given a unit vector $v \in T_xM \cap F^{-1}(1)$,
we extend it to a non-vanishing $\cC^{\infty}$-vector field $V$
on a neighborhood of $x$ in such a way that every integral curve of $V$ is geodesic,
and consider the Riemannian structure $g_V$ induced from \eqref{eq:gv}.
Then the Ricci curvature $\Ric(v)$ of $v$ with respect to $F$ coincides with the Ricci curvature
of $v$ with respect to $g_V$ (in particular, it is independent of the choice of $V$).

Let us fix a positive $\cC^{\infty}$-measure $m$ on $M$.
Inspired by the above interpretation of the Finsler Ricci curvature
and the theory of weighted Riemannian manifolds,
the weighted Ricci curvature for the triple $(M,F,m)$ was introduced in \cite{Oint} as follows.

\begin{definition}[Weighted Ricci curvature]\label{df:wRic}
Given a unit vector $v \in T_xM \cap F^{-1}(1)$,
let $\eta:(-\ve,\ve) \lra M$ be the geodesic such that $\dot{\eta}(0)=v$.
We decompose $m$ along $\eta$ using the Riemannian volume measure $\vol_{\dot{\eta}}$
of $g_{\dot{\eta}}$ as $m=e^{-\Psi}\vol_{\dot{\eta}}$, where $\Psi:(-\ve,\ve) \lra \R$.
Then we define the \emph{weighted Ricci curvature} involving a parameter $N \in [n,\infty]$ by
\begin{enumerate}[(1)]
\item $\Ric_n(v):=\displaystyle
 \begin{cases} \Ric(v)+\Psi''(0) &\ \text{if}\ \Psi'(0)=0, \\
 -\infty &\ \text{if}\ \Psi'(0) \neq 0, \end{cases}$

\item $\Ric_N(v):=\Ric(v) +\Psi''(0) -\displaystyle\frac{\Psi'(0)^2}{N-n}\quad$
for $N \in (n,\infty)$,

\item $\Ric_{\infty}(v):=\Ric(v) +\Psi''(0)$.
\end{enumerate}
We also set $\Ric_N(cv):=c^2 \Ric_N(v)$ for $c \ge 0$.
\end{definition}

We will say that $\Ric_N \ge K$ holds for some $K \in \R$ if $\Ric_N(v) \ge KF(v)^2$ for all $v \in TM$.
Observe that $\Ric_N(v) \le \Ric_{N'}(v)$ for $N<N'$,
and that for the scaled space $M'=(M,F,am)$ with $a>0$ we have $\Ric^{M'}_N(v)=\Ric^M_N(v)$.
It was shown in \cite[Theorem~1.2]{Oint} that $\Ric_N \ge K$ is equivalent to
Lott, Sturm and Villani's \emph{curvature-dimension condition} $\CD(K,N)$.
(Roughly speaking, the curvature-dimension condition is a convexity condition
of an entropy functional on the space of probability measures;
we refer to \cite{StI}, \cite{StII}, \cite{LV1}, \cite{LV2} and \cite[Part~III]{Vi}
for details and further theories.)
This equivalence extends the corresponding result on (weighted) Riemannian manifolds,
and has many analytic and geometric applications (see \cite{Oint}).

\section{Proof of Theorem~\ref{th:Funk} (Funk case)}\label{sc:Funk}

Let us first treat the Funk case.
In this section, we will denote the Funk metric simply by $F$ for brevity,
and consider the standard coordinate of $D \subset \R^n$.
The following lemma enables us to translate all the vertical derivatives ($\del/\del v^i$)
to the horizontal derivatives ($\del/\del x^i$).

\begin{lemma}{\rm (\cite[Proposition~1]{Ok}, \cite[Lemma~2.3.1]{Shspr})}\label{lm:Ok}
For any $v \in TD \setminus 0$ and $i=1,2,\ldots,n$, we have
\[ \frac{\del F}{\del x^i}(v)=F(v) \frac{\del F}{\del v^i}(v). \]
\end{lemma}

Observe that, on $TD \setminus 0$,
\begin{equation}\label{eq:FF}
\frac{1}{2} \frac{\del^2(F^2)}{\del v^i \del v^j}
 =\frac{\del}{\del v^i} \bigg[ \frac{\del F}{\del x^j} \bigg]
 =\frac{\del}{\del x^j} \bigg[ \frac{1}{F} \frac{\del F}{\del x^i} \bigg]
 = \frac{1}{F} \frac{\del^2 F}{\del x^i \del x^j}
 -\frac{1}{F^2} \frac{\del F}{\del x^i} \frac{\del F}{\del x^j}.
\end{equation}
Now, we fix a unit vector $v \in T_xD \cap F^{-1}(1)$ and choose an appropriate coordinate
that $x$ is the origin, $v=\del/\del x^n$ and that $g_{in}(v)=0$ for all $i=1,2,\ldots,n-1$.
We remark that such a coordinate exchange multiplies the Lebesgue measure
merely by a positive constant, so that the weighted Ricci curvature does not change.
Put $V:=\del/\del x^n$ on $D$ and recall that the all integral curves of $V$ are minimizing
(and hence re-parametrizations of geodesics).
Therefore it suffices to calculate the weighted Ricci curvature of $(D,g_V,m_{\bL})$.

We can represent $\del D \cap \{ x \in \R^n \,|\, x^n>0 \}$ as the graph
of the $\cC^{\infty}$-function $h:U \lra (0,\infty)$ for a sufficiently small
neighborhood $U \subset \R^{n-1}$ of $0$, namely
\begin{equation}\label{eq:h}
\del D \cap \{ (z,t) \in \R^{n-1} \times \R \,|\, z \in U,\ t>0 \}
 =\big\{ \big( z,h(z) \big) \,|\, z \in U \big\}.
\end{equation}
Then \eqref{eq:FHFF} yields
\[ F\big( V(z,t) \big) =\frac{1}{h(z)-t} \qquad \text{for}\
 (z,t) \in D \subset \R^{n-1} \times \R. \]
Putting $\del_i:=\del/\del x^i$ for simplicity, we deduce from \eqref{eq:FF} that
\begin{align*}
g_{ij}(V) &=(h-t) \del_i \del_j \bigg( \frac{1}{h-t} \bigg)
 -(h-t)^2 \del_i \bigg( \frac{1}{h-t} \bigg) \del_j \bigg( \frac{1}{h-t} \bigg) \\
&= (h-t) \bigg\{ -\frac{\del_i \del_j (h-t)}{(h-t)^2} +\frac{2\del_i (h-t) \del_j (h-t)}{(h-t)^3} \bigg\}
 -\frac{\del_i (h-t) \del_j (h-t)}{(h-t)^2} \\
&= -\frac{\del_i \del_j (h-t)}{h-t} +\frac{\del_i (h-t) \del_j (h-t)}{(h-t)^2},
\end{align*}
where the evaluations at $(z,t) \in D$ were omitted.
We remark that, for $i,j \neq n$,
\[ g_{ij}(V)=-\frac{\del_i \del_j h}{h-t}+\frac{\del_i h \del_j h}{(h-t)^2}, \quad
 g_{in}(V)=-\frac{\del_i h}{(h-t)^2}, \quad g_{nn}(V)=\frac{1}{(h-t)^2}. \]
Hence, when differentiating $g_{ij}(V(z,t))$ by $t$,
we need to take only the denominators into account.
Thus we find
\[ \frac{\del[g_{ij}(V)]}{\del t}
 =-\frac{\del_i \del_j (h-t)}{(h-t)^2} +\frac{2\del_i (h-t) \del_j (h-t)}{(h-t)^3}
 =\frac{1}{h-t}\bigg\{ g_{ij}(V) +\frac{\del_i (h-t) \del_j (h-t)}{(h-t)^2} \bigg\}. \]
Decomposing $m_{\bL}$ as $m_{\bL}=e^{-\Psi} \sqrt{\det(g_{ij}(V))}\, dx^1 dx^2 \cdots dx^n$
along the curve $\eta(t)=(0,t) \in D$, we observe
\[ \Psi(t)=\frac{1}{2} \log\Big(\! \det \big( g_{ij}(t) \big) \Big), \qquad
 \Psi'(t)=\frac{1}{2} \trace\Big[ \big( g^{ij}(t) \big) \cdot \big( g'_{ij}(t) \big) \Big], \]
where we abbreviated as $g_{ij}(t):=g_{ij}(V(0,t))$
and $(g^{ij}(t))$ stands for the inverse matrix of $(g_{ij}(t))$.
Dividing $\Psi'(t)$ by the speed $F(\dot{\eta}(t))=F(V(0,t))=(h(0)-t)^{-1}$, we obtain
\[ \big( h(0)-t \big) \Psi'(t) =\frac{1}{2}\trace\bigg[ \big( g^{ij}(t) \big) \cdot
 \bigg( g_{ij}(t) +\frac{\del_i (h(0)-t) \del_j (h(0)-t)}{(h(0)-t)^2} \bigg) \bigg]
 \equiv \frac{n+1}{2}, \]
where the second equality follows from the fact that $g_{in}(t)=-\del_i h(0)/(h(0)-t)^2=0$ for $i \neq n$
guaranteed by $g_{in}(v)=0$.
Therefore we conclude, as $(D,F)$ has the constant flag curvature $-1/4$,
\[ \Ric_{\infty}(v) =-\frac{n-1}{4}, \qquad
 \Ric_N(v) =-\frac{n-1}{4}-\frac{(n+1)^2}{4(N-n)}. \]
$\qedd$

\section{Proof of Theorem~\ref{th:Hilb} (Hilbert case)}\label{sc:Hilb}

We next consider the Hilbert case, where the calculation is similar but more involved.
We will denote the Hilbert metric of $D$ by $F$ in this section.

Given a unit vector $v \in T_xD \cap F^{-1}(1)$, similarly to the previous section,
we choose a coordinate such that $x$ is the origin, $v=\del/\del x^n$
and that $g_{in}(v)=0$ for all $i=1,2,\ldots,n-1$.
Put $V:=\del/\del x^n$ again.
In addition to $h:U \lra (0,\infty)$ as in \eqref{eq:h},
we introduce the function $b:U \lra (-\infty,0)$ such that
\[ \del D \cap \{ (z,t) \in \R^{n-1} \times \R \,|\, z \in U,\ t<0 \}
 =\big\{ \big( z,b(z) \big) \,|\, z \in U \big\}. \]
Using the Funk metric $F_+$ of $D$ and its reverse $F_-(v):=F_+(-v)$,
we can write $F(V)$ as (recall \eqref{eq:FHFF})
\[ F\big( V(z,t) \big) =\frac{F_+(V(z,t)) +F_-(V(z,t))}{2}
 =\frac{1}{2}\bigg\{ \frac{1}{h(z)-t}+\frac{1}{t-b(z)} \bigg\}. \]
It follows from Lemma~\ref{lm:Ok} and $F_-(v)=F_+(-v)$ that
\[ \frac{\del F_-}{\del x^i} =-F_- \frac{\del F_-}{\del v^i}. \]
This yields that, by putting $\del_i:=\del/\del x^i$,
\begin{align*}
2\frac{\del^2(F^2)}{\del v^i \del v^j}
&=\frac{1}{2} \frac{\del^2}{\del v^i \del v^j}(F_+^2 +2F_+ F_- +F_-^2) \\
&= \frac{1}{2} \frac{\del^2(F_+^2)}{\del v^i \del v^j}
 +\frac{1}{2} \frac{\del^2(F_-^2)}{\del v^i \del v^j}
 -\frac{\del_i F_+}{F_+} \frac{\del_j F_-}{F_-} -\frac{\del_j F_+}{F_+} \frac{\del_i F_-}{F_-} \\
&\quad +\bigg( \frac{\del_i \del_j F_+}{F_+^2} -\frac{2\del_i F_+ \del_j F_+}{F_+^3} \bigg) F_-
 +\bigg( \frac{\del_i \del_j F_-}{F_-^2} -\frac{2\del_i F_- \del_j F_-}{F_-^3} \bigg) F_+.
\end{align*}
By \eqref{eq:FF} we have, omitting the evaluations at $(z,t) \in D$,
\begin{align*}
4g_{ij}(V) &= -\frac{\del_i \del_j (h-t)}{h-t} +\frac{\del_i(h-t) \del_j(h-t)}{(h-t)^2}
 -\frac{\del_i \del_j (t-b)}{t-b} +\frac{\del_i(t-b) \del_j(t-b)}{(t-b)^2} \\
&\quad -\bigg\{ \frac{\del_i(h-t)}{h-t} \frac{\del_j(t-b)}{t-b} +\frac{\del_j(h-t)}{h-t} \frac{\del_i(t-b)}{t-b} \bigg\}
 -\frac{\del_i \del_j (h-t)}{t-b} -\frac{\del_i \del_j (t-b)}{h-t} \\
&= -\{ \del_i \del_j (h-t) +\del_i \del_j (t-b) \} \bigg( \frac{1}{h-t}+\frac{1}{t-b} \bigg) \\
&\quad +\bigg\{ \frac{\del_i(h-t)}{h-t} -\frac{\del_i(t-b)}{t-b} \bigg\}
 \bigg\{ \frac{\del_j(h-t)}{h-t} -\frac{\del_j(t-b)}{t-b} \bigg\}.
\end{align*}
Note that the assumption $g_{in}(v)=0$ implies
\begin{equation}\label{eq:h+b}
\frac{\del_i h(0)}{h(0)} -\frac{\del_i b(0)}{b(0)}=0 \qquad
 \text{for}\ i=1,2,\ldots,n-1.
\end{equation}
We also observe for later convenience that, for $i,j \neq n$,
\[ 4g_{ij}(v)=-\{ \del_i \del_j h(0) -\del_i \del_j b(0) \} \bigg( \frac{1}{h(0)}-\frac{1}{b(0)} \bigg), \quad
4g_{nn}(v)=\bigg( \frac{1}{h(0)} -\frac{1}{b(0)} \bigg)^2. \]

By the same reasoning as the Funk case, the numerators can be neglected
when one differentiates $g_{ij}(V)$ by $t$.
Thus we find
\begin{align*}
4\frac{\del[g_{ij}(V)]}{\del t}
&=-\{ \del_i \del_j (h-t) +\del_i \del_j (t-b) \} \bigg\{ \frac{1}{(h-t)^2}-\frac{1}{(t-b)^2} \bigg\} \\
&\quad +\bigg\{ \frac{\del_i(h-t)}{(h-t)^2} +\frac{\del_i(t-b)}{(t-b)^2} \bigg\}
 \bigg\{ \frac{\del_j(h-t)}{h-t} -\frac{\del_j(t-b)}{t-b} \bigg\} \\
&\quad +\bigg\{ \frac{\del_i(h-t)}{h-t} -\frac{\del_i(t-b)}{t-b} \bigg\}
 \bigg\{ \frac{\del_j(h-t)}{(h-t)^2} +\frac{\del_j(t-b)}{(t-b)^2} \bigg\}.
\end{align*}
We further calculate
\begin{align*}
4\frac{\del^2 [g_{ij}(V)]}{\del t^2}
&=-\{ \del_i \del_j (h-t) +\del_i \del_j (t-b) \} \bigg\{ \frac{2}{(h-t)^3}+\frac{2}{(t-b)^3} \bigg\} \\
&\quad +\bigg\{ \frac{2\del_i(h-t)}{(h-t)^3} -\frac{2\del_i(t-b)}{(t-b)^3} \bigg\}
 \bigg\{ \frac{\del_j(h-t)}{h-t} -\frac{\del_j(t-b)}{t-b} \bigg\} \\
&\quad +\bigg\{ \frac{\del_i(h-t)}{h-t} -\frac{\del_i(t-b)}{t-b} \bigg\}
 \bigg\{ \frac{2\del_j(h-t)}{(h-t)^3} -\frac{2\del_j(t-b)}{(t-b)^3} \bigg\} \\
&\quad +2\bigg\{ \frac{\del_i(h-t)}{(h-t)^2} +\frac{\del_i(t-b)}{(t-b)^2} \bigg\}
 \bigg\{ \frac{\del_j(h-t)}{(h-t)^2} +\frac{\del_j(t-b)}{(t-b)^2} \bigg\}.
\end{align*}
We abbreviate as $g_{ij}(t):=g_{ij}(V(0,t))$ and deduce from \eqref{eq:h+b} that, for $i,j \neq n$,
\begin{align*}
4g'_{ij}(0) &= 4g_{ij}(0) \bigg( \frac{1}{h(0)}+\frac{1}{b(0)} \bigg), \quad
4g'_{in}(0) =-\bigg( \frac{\del_i h(0)}{h(0)^2} -\frac{\del_i b(0)}{b(0)^2} \bigg)
 \bigg( \frac{1}{h(0)} -\frac{1}{b(0)} \bigg), \\
4g'_{nn}(0) &= 8g_{nn}(0) \bigg( \frac{1}{h(0)}+\frac{1}{b(0)} \bigg).
\end{align*}
We also obtain
\begin{align*}
4g''_{ij}(0) &= 8g_{ij}(0) \bigg( \frac{1}{h(0)^2} +\frac{1}{h(0)b(0)} +\frac{1}{b(0)^2} \bigg) \\
&\quad +2\bigg( \frac{\del_i h(0)}{h(0)^2} -\frac{\del_i b(0)}{b(0)^2} \bigg)
 \bigg( \frac{\del_j h(0)}{h(0)^2} -\frac{\del_j b(0)}{b(0)^2} \bigg), \\
4g''_{nn}(0) &= 8g_{nn}(0) \bigg\{ 2\bigg( \frac{1}{h(0)^2} +\frac{1}{h(0)b(0)} +\frac{1}{b(0)^2} \bigg)
 +\bigg( \frac{1}{h(0)}+\frac{1}{b(0)} \bigg)^2 \bigg\}.
\end{align*}

Put $\Psi(t)=2^{-1} \log(\det(g_{ij}(t)))$ and observe
\begin{align*}
\Psi'(t) &= \frac{1}{2} \trace\Big[ \big( g^{ij}(t) \big) \cdot \big( g'_{ij}(t) \big) \Big], \\
\Psi''(t) &= \frac{1}{2} \trace\Big[ \big( g^{ij}(t) \big) \cdot \big( g''_{ij}(t) \big)
 -\big\{ \big( g^{ij}(t) \big) \cdot \big( g'_{ij}(t) \big) \big\}^2 \Big].
\end{align*}
Comparing $g_{ij}(0)$ and $g'_{ij}(0)$, we have
\[ \Psi'(0)= \frac{1}{2}
 \bigg\{ (n-1)\bigg( \frac{1}{h(0)}+\frac{1}{b(0)} \bigg) +2\bigg( \frac{1}{h(0)}+\frac{1}{b(0)} \bigg) \bigg\}
 =\frac{n+1}{2} \bigg( \frac{1}{h(0)}+\frac{1}{b(0)} \bigg). \]
It similarly holds
\begin{align*}
\frac{1}{2} \trace\Big[ \big( g^{ij}(0) \big) \cdot \big( g''_{ij}(0) \big) \Big]
&= (n-1) \bigg( \frac{1}{h(0)^2} +\frac{1}{h(0)b(0)} +\frac{1}{b(0)^2} \bigg) \\
&\quad +\frac{1}{4} \sum_{i,j=1}^{n-1} g^{ij}(0)
 \bigg( \frac{\del_i h(0)}{h(0)^2} -\frac{\del_i b(0)}{b(0)^2} \bigg)
 \bigg( \frac{\del_j h(0)}{h(0)^2} -\frac{\del_j b(0)}{b(0)^2} \bigg) \\
&\quad +2\bigg( \frac{1}{h(0)^2} +\frac{1}{h(0)b(0)} +\frac{1}{b(0)^2} \bigg)
 +\bigg( \frac{1}{h(0)}+\frac{1}{b(0)} \bigg)^2 \\
&= (n+1) \bigg( \frac{1}{h(0)^2} +\frac{1}{h(0)b(0)} +\frac{1}{b(0)^2} \bigg)
 +\bigg( \frac{1}{h(0)}+\frac{1}{b(0)} \bigg)^2 \\
&\quad +\frac{1}{4} \sum_{i,j=1}^{n-1} g^{ij}(0)
 \bigg( \frac{\del_i h(0)}{h(0)^2} -\frac{\del_i b(0)}{b(0)^2} \bigg)
 \bigg( \frac{\del_j h(0)}{h(0)^2} -\frac{\del_j b(0)}{b(0)^2} \bigg).
\end{align*}
Combining this with
\begin{align*}
&\trace\Big[ \big\{ \big( g^{ij}(0) \big) \cdot \big( g'_{ij}(0) \big) \big\}^2 \Big] \\
&=(n-1) \bigg( \frac{1}{h(0)}+\frac{1}{b(0)} \bigg)^2 +4\bigg( \frac{1}{h(0)}+\frac{1}{b(0)} \bigg)^2 \\
&\quad +\frac{g^{nn}(0)}{8} \sum_{i,j=1}^{n-1} g^{ij}(0)
 \bigg( \frac{\del_i h(0)}{h(0)^2} -\frac{\del_i b(0)}{b(0)^2} \bigg)
 \bigg( \frac{\del_j h(0)}{h(0)^2} -\frac{\del_j b(0)}{b(0)^2} \bigg)
 \bigg( \frac{1}{h(0)} -\frac{1}{b(0)} \bigg)^2 \\
&=(n+3) \bigg( \frac{1}{h(0)}+\frac{1}{b(0)} \bigg)^2
 +\frac{1}{2} \sum_{i,j=1}^{n-1} g^{ij}(0)
 \bigg( \frac{\del_i h(0)}{h(0)^2} -\frac{\del_i b(0)}{b(0)^2} \bigg)
 \bigg( \frac{\del_j h(0)}{h(0)^2} -\frac{\del_j b(0)}{b(0)^2} \bigg),
\end{align*}
we obtain
\begin{align*}
\Psi''(0) &=(n+1) \bigg( \frac{1}{h(0)^2} +\frac{1}{h(0)b(0)} +\frac{1}{b(0)^2} \bigg)
 -\frac{n+1}{2} \bigg( \frac{1}{h(0)}+\frac{1}{b(0)} \bigg)^2 \\
&= \frac{n+1}{2} \bigg( \frac{1}{h(0)^2} +\frac{1}{b(0)^2} \bigg).
\end{align*}

Therefore we have, as $F(v)=(h(0)^{-1}-b(0)^{-1})/2=1$,
\[ \frac{d}{dt} \bigg[ \frac{\Psi'(t)}{F(V(0,t))} \bigg]_{t=0}
 =\Psi''(0) -\frac{\Psi'(0)}{2} \bigg( \frac{1}{h(0)^2} -\frac{1}{b(0)^2} \bigg)
 =-\frac{n+1}{h(0)b(0)}. \]
Since
\[ 0< -\frac{1}{h(0)b(0)} \le \frac{1}{4} \bigg( \frac{1}{h(0)}-\frac{1}{b(0)} \bigg)^2 =1, \]
this yields $\Ric_{\infty}(v) \in (-(n-1),2]$.
Moreover,
\[ \Psi'(0)^2 =\frac{(n+1)^2}{4} \bigg( \frac{1}{h(0)} +\frac{1}{b(0)} \bigg)^2
 =(n+1)^2 \bigg( 1+\frac{1}{h(0)b(0)} \bigg) \in \big[ 0,(n+1)^2 \big) \]
shows
\[ \Ric_N(v) \in \bigg(\!\! -(n-1)-\frac{(n+1)^2}{N-n},2 \bigg]. \]
$\qedd$

\section{Applications and remarks}\label{sc:appl}

As mentioned in Section~\ref{sc:prel}, $\Ric_N \ge K$ is equivalent
to the curvature-dimension condition $\CD(K,N)$.
Spaces satisfying $\CD(K,N)$ enjoy a number of properties similar to Riemannian manifolds
of $\Ric \ge K$ and $\dim \le N$.
Since $\CD(K,N)$ (between compactly supported measures) is preserved under
the \emph{pointed measured Gromov-Hausdorff convergence}
of locally compact, complete metric measure spaces (\cite[Theorem~29.25]{Vi}),
we can deal with merely bounded, convex domains $D$.

\begin{corollary}\label{cr:CD}
Let $D \subset \R^n$ be a bounded convex domain with $n \ge 2$.
Then the metric measure spaces $(D,d_{\cF},m_{\bL})$ and $(D,d_{\cH},m_{\bL})$ satisfy
$\CD(K,N)$ for $N \in (n,\infty]$ with
\[ K=-\frac{n-1}{4} -\frac{(n+1)^2}{4(N-n)}, \qquad
 K=-(n-1) -\frac{(n+1)^2}{N-n}, \]
respectively, where we read $K=-(n-1)/4$ and $K=-(n-1)$ when $N=\infty$.
In particular, they satisfy
\begin{itemize}
\item[$\bullet$]
the Brunn-Minkowski inequality by $\CD(K,N)$ with $N \in (n,\infty];$
\item[$\bullet$]
the Bishop-Gromov volume comparison theorem by $\CD(K,N)$ with $N \in (n,\infty)$.
\end{itemize}
\end{corollary}

See \cite[Proposition~2.1, Theorem~2.3]{StII}
(and \cite[Theorem~30.7]{Vi}, \cite[Theorem~6.1]{OGre} as well for the case of $N=\infty$)
for the precise statements of the Brunn-Minkowski inequality and the Bishop-Gromov volume comparison.
Beyond the general theory of the curvature-dimension condition,
the weighted Ricci curvature bound implies the following.

\begin{corollary}\label{cr:BW}
Let $D \subset \R^n$ with $n \ge 2$ be a strongly convex domain such that $\del D$ is smooth.
For $K$ as in Corollary~$\ref{cr:CD}$, $(D,F_{\cF},m_{\bL})$ and $(D,F_{\cH},m_{\bL})$ satisfy
\begin{itemize}
\item[$\bullet$]
the Laplacian comparison theorem for $N \in (n,\infty);$
\item[$\bullet$]
the Bochner-Weitzenb\"ock inequality for $N \in (n,\infty]$.
\end{itemize}
\end{corollary}

See \cite[Theorem~5.2]{OShf} for the Laplacian comparison,
and \cite[Theorems~3.3, 3.6]{OSbw} for the Bochner-Weitzenb\"ock formula
(by the Bochner-Weitzenb\"ock inequality we meant the inequality given by
plugging the weighted Ricci curvature bound into the Bochner-Weitzenb\"ock formula).

We conclude the article with remarks on possible improvements of the estimates
in Theorems~\ref{th:Funk}, \ref{th:Hilb}.
Our estimates on $\Ric_N$ with respect to $m_{\bL}$ are independent of the shape of $D$.
In particular, Theorem~\ref{th:Hilb} provides the same (far from optimal) estimates
even for the Klein model of the hyperbolic spaces.
Thus there would be a better choice of a measure depending on the shape of $D$.
Then, as an arbitrary measure is represented by $e^{-\psi} m_{\bL}$,
its weighted Ricci curvature is calculated by combining Theorems~\ref{th:Funk}, \ref{th:Hilb}
and the convexity of $\psi$.
One may think of the squared distance function from some point as a candidate of $\psi$,
however, in order to estimate its convexity along geodesics,
we need to bound not only the flag curvature but also the \emph{uniform convexity}
as well as the \emph{tangent curvature} (see \cite[Theorem~5.1]{Ouni}).
The uniform convexity is measured by the constant
\[ \bC =\sup_{x \in M} \sup_{v,w \in T_xM \setminus 0} \frac{F(w)}{g_v(w,w)^{1/2}}, \]
and it is infinite for Funk metrics.
As for Hilbert geometry, one could bound $\bC$ by the convexity of $\del D$
(whereas it seems unclear; see \cite[Remark~2.1]{Eg}).
The author has no idea about the tangent curvature,
which measures how the tangent spaces are distorted as one moves in $M$.

There are several natural constructive measures $m$ on $D$,
and it is interesting to consider the corresponding weighted Ricci curvature $\Ric_N^m(V)$.
Then, however, it seems not easy (at least more difficult than $m_{\bL}$) to calculate $\Ric_N^m(V)$
because $m$ should depend on the shape of whole $\del D$, while $g_V$ is induced
only from the behavior of $F_{\cF}$ or $F_{\cH}$ near the direction $V$.

We also remark that, in Hilbert geometry (which is both forward and backward complete),
$\Ric_N$ with $N<\infty$ can not be nonnegative for any measure.
Otherwise, $g_V$ splits isometrically that is a contradiction (\cite[Proposition~4.3]{Ospl}).
Due to the same reasoning, $\Ric_{\infty}$ can be nonnegative only when $\sup \Psi=\infty$.

{\small

}

\end{document}